\documentclass[review,authoryear,12pt]{elsarticle}

\usepackage{amssymb}

\usepackage{amsmath}
\usepackage{graphicx}
\usepackage{natbib,amssymb}

\usepackage{soul}
\usepackage{xcolor}

\newtheorem{theorem}{Theorem}[section]
\newtheorem{lemma}[theorem]{Lemma}

\journal{Statistics and Probability Letters}

\begin{document}

\begin{frontmatter}

\title{Fractional Poisson processes and their representation by infinite systems of ordinary differential equations\tnoteref{t1}}
\tnotetext[t1]{The University of Adelaide preprint number for this article is: ADP-13-16/T836}

\author[CS]{Markus Kreer}
\ead{mkreer@web.de}
\address[CS]{CAMPUSERVICE GmbH, Servicegesellschaft der Johann Wolfgang Goethe-Universit{\"a}t Frankfurt, Rossertstrasse 2,60323 Frankfurt am Main, Germany}
\author[adel]{Ay{\c s}e K{\i}z{\i}lers{\" u} \corref{cor1}}
\cortext[cor1]{Corresponding author}
\ead{ayse.kizilersu@adelaide.edu.au}
\address[adel]{CSSM,
   The School of Chemistry and Physics
   Department of Physics and Mathematical Physics,
   Adelaide University, 5005, Australia}
\author[adel]{Anthony W. Thomas}

\begin{abstract}
Fractional Poisson processes, a rapidly growing area of non-Markovian stochastic processes, are useful in statistics to describe data from counting processes when waiting times are not exponentially distributed. We show that the fractional Kolmogorov-Feller equations for the probabilities at time $t$ can be representated by an infinite linear system of ordinary differential equations of first order in a transformed time variable. These new equations resemble  a linear version of the discrete coagulation-fragmentation equations, well-known from the non-equilibrium theory of gelation, cluster-dynamics and phase transitions in physics and chemistry.

\end{abstract}


\end{frontmatter}



\section{Introducing fractional Poisson processes}
Since the late 1990s there has been a great interest in non-Markovian continous-time processes, especially those arising from waiting times 
between two events  that are not exponentially distributed (for a general overview see \cite{Embrechts97}, \cite{Grandell97}) but sub-exponentially, for example 
 (\cite{Jumarie01}, \cite{Laskin03})
\begin{equation}
\text{prob}(T_w<t) = 1 - E_{\beta}(-\lambda t^{ \beta })
\label{eqn:ml_wt}
\end{equation} 
where  $\text{prob}(T_w<t)$ is the probability that the waiting time $T_w$ is less than some $t$, $\lambda>0$,  $0<\beta\leq 1$ and $E_\beta$ is Mittag-Leffler function. The Mittag-Leffler function in Eqn.(\ref{eqn:ml_wt}), has the series representation 
$E_{\beta}(z) = \sum_{m=0}^{\infty}z^{m}/\Gamma(\beta m +1) $
and is a fractional generalization of the exponential function; where $\Gamma(x)=\int_{0}^{\infty} s^{x-1} \exp{(-s)} ds$ is the Gamma function (for $\beta=1$ the exponential function is recovered).
 
Starting from the waiting time distribution Eqn.(\ref{eqn:ml_wt}),  \cite{Laskin03} introduced the fractional Poisson process 
as the counting process with probability $P_{\beta}(n,t)$ of  $n$ items ($n = 0, 1, 2, . . .$) arriving by a time $t$.  
\cite{Beghin09} pointed out, that the approach of  \cite{Laskin03} is equivalent to solving their fractional equations  
\begin{equation}
\frac{d^{\beta} }{dt^{\beta}} P_{\beta}(n,t) = \lambda (  P_{\beta}(n-1,t)  -  P_{\beta}(n,t)  )   \hspace{0.2cm} , \hspace{0.5cm} 0<\beta\leq 1 \, ,
\label{eqn:fode1}
\end{equation}
where the fractional derivative is taken in the sense of Dzerbayshan-Caputo (e.g. \cite{Podlubny99}  p. 78), and is defined on a twice continously differentiable function $f(t)$ as the usual derivative for $\beta=1$ and for
 $0<\beta<1$ is
\begin{equation}
\frac{d^{\beta} }{dt^{\beta}} f(t) \equiv \frac{1}{\Gamma(1-\beta)}    \int_{0}^{t} \frac{ds}{( t- s)^{\beta}} \frac{d }{ds} f(s)\,.
\label{eqn:Caputo}
\end{equation}
Their solution with the initial condition $P_{\beta}(n,t=0)  = \delta_{n,0}$ is given by
\begin{eqnarray}
P_{\beta}(n,t) = (-1)^{n}\sum_{j=n}^{\infty} \binom{j}{n}    \frac{(-1)^{j} (\lambda t^{\beta})^{j}  }{ \Gamma( \beta j +1)} = (-1)^{n}\sum_{j=n}^{\infty} \binom{j}{n}  \phi_{\beta}(j,t)\, 
\label{eqn:p_solution}
\end{eqnarray}
where the functions $\phi_{\beta}(j,t)$ are defined as
\begin{equation}
\phi_{\beta}(j,t) \equiv (-1)^{j}\frac{ (\lambda t^{\beta})^{j}}{ \Gamma( \beta j +1)} \hspace{0.2cm} ,  \hspace{0.5cm}   j=0,1,2,...
\label{eqn:phi_def}
\end{equation}
Note that for $n=0$  in Eqn.(\ref{eqn:p_solution}) we recover the Mittag-Leffler function with index $\beta$, $ E_{\beta}(-\lambda t^{ \beta })$. This is not surprising as
$P_{\beta}(0,t)$ is the probability of no birth taking place up to time $t$, i.e. $\text{prob}(T_w>t)$ which is nothing else but the complement to the waiting time probability, Eqn.(\ref{eqn:ml_wt}).
Furthermore, one readily recovers the normalization condition $\sum_{n=0}^{\infty} P_{\beta}(n,t) =1$.

Equation (\ref{eqn:fode1}) together with Eqn.(\ref{eqn:Caputo}) are now considered the standard fractional Kolmogorov-Feller equations for a fractional Poisson process. They are the basis of a fast growing branch of probability theory (e.g. \cite{Beghin09}, \cite{Beghin10}, \cite{Beghin13}, \cite{Orsingher13}). 

We shall show in this letter that the fractional Poisson process Eqn.(\ref{eqn:fode1}) can also be described by an infinite linear system of ordinary differential equations (ODE) of first order in a transformed time variable for the probabilities $P_{\beta}(n,t)$  on the left-hand side, and on the right-hand side we have the usual two terms of a standard Poisson process plus infinitely many more terms consisting of some time-independent constants times $P_{\beta}(m,t)$,  and  $m>n$. This result, formulated in Theorem 3.1, Eqn.(\ref{eqn:birth-death}), is what we term the ``ODE-representation of the Kolmogorov-Feller equations''. 

These new equations for the fractional Poisson process bear a striking resemblance to the linear version of the discrete cluster equations for the Glauber kinetic Ising model, as discussed in \cite{Binder74} and \cite{Kreer93} which is  briefly discussed in Appendix B,  Eqn.(\ref{eqn:lin_coag_frag}).
These cluster ODEs describe a typical dynamics in which clusters consisting of $n$ particles, say, can coagulate with other clusters to form larger clusters or fragment to form smaller ones. Thus, our approach allows us, in principle, to understand the dynamics of the fractional Poisson process in terms of cluster interactions.
Consequently, our ``ODE-representation of the Kolmogorov-Feller equations"  belongs to a wider class of coagulation-fragmentation equations, which are of major interest to the wider community,
as they show a variety of interesting features, such as asymptotic self-similarity for large times (``dynamical scaling''), gelating at finite time, metastablility, etc. They may also be a tool to understand phenomena of non-equilibrium statistical physics and phase transitions.  For further reading we refer to  \cite{Spouge84}, \cite{BallCarr90}, \cite{Kreer93}, \cite{daCosta95}, \cite{Mischler02} and \cite{Lamb10}. 

One of the aims of this letter is to bring the community of probabilists dealing with fractional Poisson processes and the community of analysts and physicists dealing with coagulation-fragmentation equations closer as both subjects seem to be more related than previous thought.
 
The plan of the paper is as follows: Before presenting our main result in Section 3 we  prove in Section 2 some necessary Lemma dealing with a combinatorial inversion formula for certain functions.

\section{An application of Krattenthalers theorem}
To prove an essential Lemma, we 
heuristically define ``infinite'' vectors from the solutions Eqns.(\ref{eqn:p_solution}) and~(\ref{eqn:phi_def}) of the fractional Poisson process,
\begin{eqnarray}
               & &[P_{\beta}(0,t),-P_{\beta}(1,t),\cdots, (-1)^{n} P_{\beta}(n,t),\cdots]^T     \nonumber \\
\text{and}&& [\phi_{\beta}(0,t),\phi_{\beta}(1,t),\cdots,  \phi_{\beta}(n,t),\cdots]^T . \nonumber
\end{eqnarray}
We  may then write the solutions Eqn.(\ref{eqn:p_solution}) formally  as follows
\begin{eqnarray}
\left[ \begin{array}{c}
~P_{\beta}(0,t) \\ -P_{\beta}(1,t)  \\ ~P_{\beta}(2,t) \\ -P_{\beta}(3,t) \\ ~P_{\beta}(4,t)  \\ \vdots  
\end{array} \right]
 = \left[ \begin{array}{cccccc}
1 & 1 & 1 & 1 & 1 &  \cdots \\
0 & 1 & 2 & 3 & 4 &  \cdots \\
0 & 0 & 1 & 3 & 6 &  \cdots \\
0 & 0 & 0 & 1 & 4 &  \cdots \\
0 & 0 & 0 & 0 & 1 &  \cdots \\
\vdots & \vdots &\vdots & \vdots &\vdots & \ddots 
\end{array} \right]
\left[ \begin{array}{c}
\phi_{\beta}(0,t) \\  \phi_{\beta}(1,t)  \\ \phi_{\beta}(2,t) \\ \phi_{\beta}(3,t) \\  \phi_{\beta}(4,t)  \\ \vdots  
\end{array} \right] \, ,
\label{eqn:matrixeqn}
\end{eqnarray}
where the matrix is infinite, triangular and contains the binomial coefficients as non-zero elements. 
This matrix in Eqn.(\ref{eqn:matrixeqn}) is well known as Pascal matrix and its inverse is readily known from \cite{Rita02} as
\begin{eqnarray}
 \left[ \begin{array}{rrrrrr}
1 & -1 & 1 & -1 & 1 &  \cdots \\
0 & 1 & -2 & 3 & -4 &  \cdots \\
0 & 0 & 1 & -3 & 6 &  \cdots \\
0 & 0 & 0 & 1 & -4 &  \cdots \\
0 & 0 & 0 & 0 & 1 &  \cdots \\
\vdots & \vdots &\vdots & \vdots &\vdots & \ddots 
\end{array} \right] \, . \nonumber
\end{eqnarray}
Thus formal inversion of Eqn.(\ref{eqn:matrixeqn}) yields 
\begin{eqnarray}
\left[ \begin{array}{c}
\phi_{\beta}(0,t) \\  \phi_{\beta}(1,t)  \\ \phi_{\beta}(2,t) \\ \phi_{\beta}(3,t) \\  \phi_{\beta}(4,t)  \\ \vdots  
\end{array} \right] =
\left[ \begin{array}{rrrrrr}
1 & -1 & 1 & -1 & 1 &  \cdots \\
0 & 1 & -2 & 3 & -4 &  \cdots \\
0 & 0 & 1 & -3 & 6 &  \cdots \\
0 & 0 & 0 & 1 & -4 &  \cdots \\
0 & 0 & 0 & 0 & 1 &  \ddots \\
\vdots & \vdots &\vdots & \vdots &\vdots & \ddots 
\end{array} \right]
\left[ \begin{array}{c}
~P_{\beta}(0,t) \\ -P_{\beta}(1,t)  \\ ~P_{\beta}(2,t) \\ -P_{\beta}(3,t) \\ ~P_{\beta}(4,t)  \\ \vdots  
\end{array} \right]\, ,
\nonumber
\end{eqnarray}
or in component form this is
\begin{eqnarray}
\phi_{\beta}(m,t) = (-1)^{m}\sum_{k=m}^{\infty} \binom{k}{m}  P_{\beta}(k,t)\,.
\label{eqn:p_solutionII}
\end{eqnarray}
For $m=0$ we recover the normalization condition of the probabilities. For $m>0$ the right-hand side of equation Eqn.(\ref{eqn:p_solutionII}) 
is bounded  by the factorial moments, which are known to be finite and given by \cite{Beghin10}, Eqn.(2.32).

For the rigorous path to obtain Eqn.(\ref{eqn:p_solutionII}), we may just apply the Krattenthaler theorem from \cite{Krattenthaler96}.
This rather useful result from combinatorics was one of the first in the theory of Riordan arrays and handles the inversion of certain classes of infinite matrices and related ``rotated inversion formula''. Following the original notation in \cite{Krattenthaler96}, Eqns.(1.1)(1)~and ~(1.1)(2), we chose
$a_{j} = (j+1)$ and $b_{j}=0$, and recover both our binomial coefficients in the original infinite matrix as well as the binomial 
coefficients with alternating signs in its inverse matrix. Eqn.(4.2) in \cite{Krattenthaler96} now provides the desired  ``rotated inversion formula'' (originally due to \cite{Riordan68} if inverse of matrix is already known) and corresponds to Eqn.(\ref{eqn:p_solutionII}). Note that all our infinite series converge (see \cite{Beghin10}) and we have thus proven the following
\begin{lemma}
\emph{(Rotated Inversion Formula)}
\label{thm:rotatedinversion}
For $\beta \in (0,1]$ and $j=0,1,2,...$ define the family of real-valued functions in $t$, $\phi_{\beta}(j,t)= (-\lambda t^{\beta})^{j}/\Gamma(\beta j +1)$ and let $P_{\beta}(n,t)$ for $n=0,1,2,...$ denote the solutions of  Eqn.(\ref{eqn:fode1}) with initial conditions $P_{\beta}(n,t=0)  = \delta_{n,0}$ . Then the following two identities  imply each other
\begin{eqnarray}
P_{\beta}(n,t) &=&(-1)^{n}\sum_{j=n}^{\infty} \binom{j}{n}  \phi_{\beta}(j,t) \, ,\nonumber \\
\phi_{\beta}(m,t) &=& (-1)^{m}\sum_{k=m}^{\infty} \binom{k}{m}  P_{\beta}(k,t)\, .
\nonumber
\end{eqnarray}
\end{lemma}

\section{Main result}

\begin{theorem}
\emph{(ODE-representation of Kolmogorov-Feller equations)}
\label{thm:ODE-representation}
After the transformation $\tau \equiv t^{\beta}$, the fractional equations  Eqn.(\ref{eqn:fode1}) with solution $P_{\beta}(n,t)$, $n=0,1,2,...$ defined in  Eqn.(\ref{eqn:p_solution}) are equivalent  
to the following infinite system of ordinary differential equations
\begin{equation}
\frac{\partial }{\partial \tau} P_{\beta}(n,\tau) = \sum_{k=n-1}^{\infty} A_{n,k}  P_{\beta}(k,\tau)   \hspace{0.2cm},\hspace{0.5cm} n=0,1,2,...
\label{eqn:birth-death}
\end{equation}
whereby the constant coefficients vanish for $k<n-1$ by the properties of binomial coefficients and for $k\geq n-1$  are defined by
\begin{eqnarray}
 A_{n,k} &= & (-1)^{n+1} \lambda \sum_{j=n-1}^{k}  (-1)^{j}   \binom{j+1}{n} \binom{k}{j} \frac{\Gamma( \beta j+1 )}{\Gamma( \beta j +\beta)}  
\label{eqn:coefficients}
\end{eqnarray}
By construction, the functions $P_{\beta}(n,\tau)$, $n=0,1,2,...$ defined in Eqn.(\ref{eqn:p_solution}) with $\tau\equiv t^{\beta}$ solve 
Eqn.(\ref{eqn:birth-death}) with initial conditions $P_{\beta}(n,\tau=0)  = \delta_{n,0}$.
\end{theorem}

{\bf Proof.}
Differentiating for each $n=0,1,2,...$ the functions $P_{\beta}(n,t)$  in Eqn.(\ref{eqn:p_solution}) with respect to $t$ and remembering that a convergent power series is differentiable in the domain of convergence componentwise, we obtain
\begin{eqnarray}
\frac{\partial }{\partial t} P_{\beta}(n,t) = (-1)^{n}\sum_{j=n}^{\infty} \binom{j}{n}    \frac{(-1)^{j} (\lambda t^{\beta})^{j-1}  }{ \Gamma( \beta j +1)} j \lambda \beta t^{\beta -1}\, .
\label{eqn:inf_aux}
\end{eqnarray}
Next dividing Eqn.(\ref{eqn:inf_aux}) by $\beta t^{\beta-1}$ and defining $\tau \equiv t^{\beta}$ we get
\begin{subequations}
\begin{eqnarray}
\frac{\partial }{\partial \tau} P_{\beta}(n,\tau) &=& -(-1)^{n}\lambda\sum_{j=n}^{\infty} \binom{j}{n}    \frac{(-1)^{j-1} (\lambda \tau)^{j-1}  }{ \Gamma( \beta j )}   \label{eqn:line1} \\
& = & (-1)^{n+1}\lambda\sum_{j=n}^{\infty} \binom{j}{n} \frac{\Gamma( \beta (j-1)+1 )}{\Gamma( \beta j )}   \frac{(-1)^{j-1} (\lambda \tau)^{j-1}  }{ \Gamma( \beta (j-1)+1 )}  \\
& = & (-1)^{n+1}\lambda\sum_{j=n}^{\infty} \binom{j}{n} \frac{\Gamma( \beta (j-1)+1 )}{\Gamma( \beta j )}   
 \phi_{\beta}(j-1,\tau)  \\
& = & (-1)^{n+1}\lambda\sum_{j=n-1}^{\infty} \binom{j+1}{n} \frac{\Gamma( \beta j+1 )}{\Gamma( \beta j +\beta)}   
 \phi_{\beta}(j,\tau) \\
& = & (-1)^{n+1}\lambda\sum_{j=n-1}^{\infty} \binom{j+1}{n} \frac{\Gamma( \beta j+1 )}{\Gamma( \beta j +\beta)}   \label{eq:line5}
    \sum_{k=j}^{\infty} \binom{k}{j}  (-1)^{j}    P_{\beta}(k,\tau)                 \nonumber \\
    & & \\
& = & (-1)^{n+1}\lambda\sum_{k=n-1}^{\infty}  P_{\beta}(k,\tau)  \sum_{j=n-1}^{k}  (-1)^{j}   \binom{j+1}{n} \binom{k}{j} \frac{\Gamma( \beta j+1 )}{\Gamma( \beta j +\beta)} \, ,    \nonumber \\
\label{eqn:inf_odes}
\end{eqnarray}
\end{subequations}
where we have used the convenient notation $ P_{\beta}(-1,\tau) =0$ and applied our {rotated inversion Lemma}~\ref{thm:rotatedinversion} in Eqn.\eqref{eq:line5}. We also interchanged the order of summation in the double sum in the last line, Eqn.\eqref{eqn:inf_odes} using Cauchy's double series theorem (see Appendix A for the necessary proof of absolute convergence).  $\blacksquare$

Note that our change of variable  $\tau \equiv t^{\beta}$, $0<\beta\leq 1$ in equation Eqn.(\ref{eqn:line1}) is motivated by convenience and arises in a natural way. From Eqn.~(\ref{eqn:p_solution}) we see immediately that $P_{\beta}(n,\tau)$ is at least once continuosly differentiable with respect to $\tau$ because the factorial moments are known to be finite (\cite{Beghin10}, Eqn.(2.32)). Thus we can differentiate using the chain rule to obtain
\begin{eqnarray}
\frac{\partial }{\partial t} P_{\beta}(n,\tau=t^{\beta})=  \frac{\partial }{\partial \tau} P_{\beta}(n,\tau) \cdot  \frac{\partial }{\partial t}t^{\beta}\,.
\label{eqn:my_way}
\end{eqnarray}
By contrast, \cite{Jumarie01} requires in his ``non-standard analysis of order $\beta$" the $\beta$-continuity for a real valued function $y(t)$ (his definition 2.1), $y(t+\theta)-y(t)  = {o}(\theta^{\beta})$, and the $\beta$ finite derivative to be $\left(d/dt^{\beta}\right) y(t)=\lim_{\theta \downarrow 0} \left[\theta^{-\beta}(y(t+\theta)-y(t)) \right] $ (his definition 2.2). Applying both definitions to the probability $P_{\beta}(n,\tau)$ in his Section 5.,
he obtains a different solution for $P_{\beta}(n,\tau)$ because his approach supreses the factor $(\partial / \partial t) t^{\beta}$ of our equation  Eqn.(\ref{eqn:my_way}).  In his article \cite{Jumarie01} points out that his non-standard analysis is different from the analysis using fractional derivatives and consequently leading to different results, namely a solution similar to the standard Poisson process albeit with a stretched exponential.

\section{Acknowledgement}
This work was supported by the University of Adelaide and the Australian Research Council (FL 0993347). 

\begin{appendix}
\section{ Absolute convergence of double series}
Defining the generating function for $|s|<1$,
\begin{equation}
 G_{\beta}(s,t) = \sum_{n=0}^{\infty} s^{n} P_{\beta}(n,t)\, ,
\label{eqn:genfunc1}
\end{equation}
we note that this is a non-negative function for $0\leq s < 1$. It has been derived in \cite{Beghin10}, Section 2.3  (and also before by \cite{Laskin03} for his version of the fractional process) and is in fact the celebrated Mittag-Leffler function with modified arguments, namely
\begin{equation}
 G_{\beta}(s,t) = E_{\beta}(\lambda t^{\beta}(s-1) )\, .
\label{eqn:solution}
\end{equation}
From the generating function Eqn.(\ref{eqn:genfunc1}) and its explicite solution Eqn.(\ref{eqn:solution}) we can chose a constant $\sigma>e$  and a positive function $ F(\tau,\sigma)$  such that $ \left|  P_{\beta}(k,\tau) \right| < e^{-k \sigma} F(\tau,\sigma)$ for any $k=0,1,2,...$. Thus equipped we start to show the absolute convergence for 
any fixed integer $n\geq 0$
\begin{eqnarray}
& & \sum_{j=n-1}^{\infty} \left|   \binom{j+1}{n} \frac{\Gamma( \beta j+1 )}{\Gamma( \beta j +\beta)}   
    \sum_{k=j}^{\infty} \binom{k}{j}  (-1)^{j}    P_{\beta}(k,\tau) \right|                \nonumber \\
& \leq &
\sum_{j=n-1}^{\infty} \left|   \binom{j+1}{n} \frac{\Gamma( \beta j+1 )}{\Gamma( \beta j +\beta)}   
    \sum_{k=j}^{\infty} \binom{k}{j}    e^{-k \sigma} F(\tau,\sigma) \right|                \nonumber \\
& \leq &
 F(\tau,\sigma) \sum_{j=n-1}^{\infty} \left|   \frac{(j+1)}{n! \Gamma(j+1-n+1)} \frac{\Gamma( \beta j+1 )}{\Gamma( \beta j +\beta)}   
    \sum_{k=j}^{\infty}  \frac{\Gamma( k+1 )}{\Gamma( k+1-j)}      e^{-k \sigma} \right|    \,.            \nonumber
\label{eqn:abs_convergence}
\end{eqnarray}
With the following inequality for all $k\geq j$ for any $j\geq 1$ fixed
\begin{equation}
\frac{\Gamma( k+1 )}{\Gamma( k+1-j)} = k (k-1) (k-2) \cdots (k-j) < (k+1)^{j}\, ,
\label{eqn:trivial_quotient}
\end{equation}
we find a suitable  constant $C_{1}>0$ independent of $j$ such that by the integral test of convergence
\begin{eqnarray}
\left|   \sum_{k=j}^{\infty}  \frac{\Gamma( k+1 )}{\Gamma( k+1-j)}      e^{-k \sigma} \right|   
&< & C_{1} \sum_{k=j}^{\infty}  ( k+1 )^{j}   e^{-k \sigma}  \nonumber \\
& \leq &  \tilde{C}_{1} e^{\sigma} \int_{0}^{\infty} \frac{x^{j}}{\sigma^{j}}  e^{-x} \frac{dx}{\sigma} 
 \leq   \tilde{C}_{1} e^{\sigma} \Gamma{(j+1)} e^{- (\ln{\sigma})\cdot j} \ , \nonumber
\label{eqn:bound1}
\end{eqnarray}
where $\tilde{C}_{1}=C_1\cdot \sigma$. Finally with this result we immediately conclude that because of $\ln{\sigma}>0$
\begin{eqnarray}
& & \sum_{j=n-1}^{\infty} \left|   \binom{j+1}{n} \frac{\Gamma( \beta j+1 )}{\Gamma( \beta j +\beta)}   
    \sum_{k=j}^{\infty} \binom{k}{j}  (-1)^{j}    P_{\beta}(k,\tau) \right|                \nonumber \\
& \leq &
 F(\tau,\sigma) \tilde{C_{1}} e^{\sigma} \sum_{j=n-1}^{\infty}   \frac{\Gamma(j+2)}{n! \Gamma(j+1-n+1)} \frac{\Gamma( \beta j+1 )}{\Gamma( \beta j +\beta)}   e^{- (\ln{\sigma})\cdot j}             \nonumber \\
& \leq &
 F(\tau,\sigma) \tilde{C}_{1} K e^{\sigma} \sum_{j=0}^{\infty}  (j+2)^{n-1} (j+1) (\beta j+1)^{1- \beta}   e^{- (\ln{\sigma})\cdot j}   \hspace{0.2cm} <  \hspace{0.2cm}  \infty   \ , \nonumber 
\label{eqn:bound2}
\end{eqnarray}
where we have again applied Eqn.(\ref{eqn:trivial_quotient}) on the quotient of Gamma functions with integer argument and bounded the second quotient
of Gamma functions by $K (\beta j +1)^{1-\beta}$ with a suitable $K>0$ due to the asymptotics of Weierstrass, 
$\Gamma(x+1)/\Gamma(x+1-\epsilon) \sim x^{\epsilon}$ for $x\rightarrow \infty$ (\cite{ParisKaminski01}, Eqn.(2.2.30)). $\blacksquare$

\section{Some remarks about cluster equations and their linear version}

Let $c_n(\cdot)$  denote the concentration of $n$-clusters in a physical system of magnetic spins, say ( i.e.  elementary magnets pointing in one direction), which can form bigger clusters with the same magnetization. Thus, the concentration corresponds for small values to the probability $P(n,t)$ to find an $n$-cluster in the unit volume at some time $t>0$.
Then the cluster equations of \cite{Binder74} describing the dynamics, read in the notation of \cite{BallCarr90} (see also \cite{Kreer93}) as
\begin{eqnarray}
\dot{c}_n &=& \frac{1}{2} \sum_{k=1}^{n} W_{n-k,k-1} - \sum_{k=1}^{\infty} W_{n,k-1}
\hspace{0.2cm} , \hspace{0.2cm} n=1,2,... ,
\label{eqn:coag_frag} 
\end{eqnarray}
where $c_0(t)\equiv 1$ and the cluster currents $W_{n,k}$ are defined as 
\begin{eqnarray}
W_{n,k} & = & a_{n,k} z c_{n}c_{k} - b_{n,k}c_{n+k+1}  \,.
\label{eqn:cluster_current}
\end{eqnarray}
Note that due to physical considerations the time-independent coefficients in  Eqn.(\ref{eqn:cluster_current}) are symmetric, $a_{n,k}=a_{k,n}\geq 0$, $b_{n,k}=b_{k,n}\geq 0$ (see \cite{BallCarr90} an references therein). The constant $z>0$ corresponds to the concentration of $1$-clusters provided by an external particle reservoir. Furthermore, in the expression for the cluster current $W_{n,k}$ in Eqn.(\ref{eqn:cluster_current})
\begin{itemize}
\item the first term describes a coagulation process when a $n$-cluster coagulates with a $k$-cluster via a newly created $1$-cluster linking both of them to form a $(n+k+1)$-cluster,
\item and the second term describes a fragmentation process when a $n+k+1$-cluster fragments to an $n$- and a $k$-cluster via deletion of a $1$-cluster.
\end{itemize}
Inserting now in Eqn.(\ref{eqn:cluster_current}) the special choice for 
\begin{eqnarray}
a_{n-k,k-1} = \begin{cases}
a_{n-1,0} & \quad \textrm{if}\quad k=1 \\
0 & \quad \textrm{if} \quad k\neq1
\end{cases}
\nonumber
\end{eqnarray} 
we obtain from Eqn.(\ref{eqn:coag_frag}) the aforementioned linear version\footnote{Note that for $b_{n-k,k-1}=b_{n-1,0}$ for $k=1$ and $0$ otherwise, we recover the equations for birth-death processes.}
\begin{eqnarray}
\dot{c}_n &=& \left( a_{n-1,0}z \right) c_{n-1} - \left( a_{n,0}z + \frac{1}{2}\sum_{k=1}^{n}b_{n-k,k-1}\right) c_n  + \sum_{k=1}^{\infty}b_{n,k-1} c_{n+k} \nonumber
\\
\label{eqn:lin_coag_frag}
\end{eqnarray}

\end{appendix}




\bibliographystyle{elsarticle-harv}

\end{document}